%% file: agt-4-1.tex
\documentclass{gtart}

\input agtout

\lognumber{1}
\volumenumber{4}
\volumeyear{2004}
\papernumber{1}
\published{9 January 2004}
\pagenumbers{1}{22}
\received{27 July 2003}
\revised{3 January 2004}
\accepted{7 January 2004}

\usepackage{amssymb, amsmath, epsf}

\def\zz{{\bf Z}}
\def\qq{{\bf Q}}

\def\calc{\mathcal{C}}
\def\calg{\mathcal{G}}
\def\cala{\mathcal{A}}

\def\cfigure#1#2#3{\begin{figure}[ht!]
\cl{\epsfxsize=#1
\epsfbox{#2}}
\caption{#3}
\end{figure}}

\newtheorem{theorem}{Theorem}[section]
\newtheorem{lemma}[theorem]{Lemma}

\theoremstyle{definition}
\newtheorem{definition}[theorem]{Definition}
\newtheorem{example}[theorem]{Example}

\numberwithin{equation}{section}

\begin{document}

\title{The concordance genus of knots}
\author{Charles Livingston}
\address{Department of Mathematics, Indiana University\\Bloomington, 
IN 47405, USA}
\email{livingst@indiana.edu}

\begin{abstract}  
In knot concordance three genera arise naturally, $g(K), g_4(K)$, and
$g_c(K)$: these are the classical genus, the 4--ball genus, and the
concordance genus, defined to be the minimum genus among all knots
concordant to $K$.  Clearly $0 \le g_4(K) \le g_c(K) \le g(K)$.
Casson and Nakanishi gave examples to show that $g_4(K)$ need not
equal $g_c(K)$.  We begin by reviewing and extending their results.

For knots representing elements in $\cala$, the concordance group of
algebraically slice knots, the relationships between these genera are
less clear.  Casson and Gordon's result that $\cala$ is nontrivial
implies that $g_4(K)$ can be nonzero for knots in $\cala$.  Gilmer
proved that $g_4(K)$ can be arbitrarily large for knots in $\cala$.
We will prove that there are knots $K$ in $\cala$ with $g_4(K) = 1$
and $g_c(K)$ arbitrarily large.

Finally, we tabulate $g_c$ for all prime knots with 10 crossings and,
with two exceptions, all prime knots with fewer than 10
crossings. This requires the description of previously unnoticed
concordances.
\end{abstract}

\asciiabstract{%
In knot concordance three genera arise naturally, g(K), g_4(K), and
g_c(K): these are the classical genus, the 4-ball genus, and the
concordance genus, defined to be the minimum genus among all knots
concordant to K.  Clearly 0 <= g_4(K) <= g_c(K) <= g(K).  Casson and
Nakanishi gave examples to show that g_4(K) need not equal g_c(K).
We begin by reviewing and extending their results.

For knots representing elements in A, the concordance group of
algebraically slice knots, the relationships between these genera are
less clear.  Casson and Gordon's result that A is nontrivial implies
that g_4(K) can be nonzero for knots in A.  Gilmer proved that g_4(K)
can be arbitrarily large for knots in A.  We will prove that there are
knots K in A with g_4(K) = 1 and g_c(K) arbitrarily large.

Finally, we tabulate g_c for all prime knots with 10 crossings and,
with two exceptions, all prime knots with fewer than 10 crossings.
This requires the description of previously unnoticed concordances.}

 
\maketitle
 
\section{Introduction and basic results}

For a knot $K \subset S^3$,  three genera  arise naturally: $g(K)$,  the genus
of $K$,   is the minimum genus among surfaces bounded by $K$ in $S^3$; $g_4(K)$ is the
minimum genus among surfaces bounded by $K$ in $B^4$; $g_c(K)$ is the minimum value of
$g(K')$ among all knots $K'$ concordant to $K$.  This paper investigates the relationships
between these knot invariants.

The classical genus  came to be fairly  well understood through Schubert's work \cite{sc}
proving that knot genus is additive under connected sum.  The 4--ball genus is far more
subtle.  Even the fact that 
$g_4(K)$ can be zero for a nontrivial knot is not entirely obvious; this was seen first as a
consequence of    Artin's construction of a knotted $S^2$ in $S^4$ \cite{ar}. That
$g_4(K)$ can be nonzero for a nontrivial knot was first proved by Fox and Milnor \cite{fo,
fm} and by Murasugi,
\cite{mu}. The 4-ball genus remains an object of investigation;    for instance, the
solution of the Milnor  conjecture, proved in 
\cite{km}, implies that  for a torus knot $K$,  $g_4(K) = g(K)$,  in the smooth
category.  (This is false in the topological locally flat category, as observed by Rudolph~\cite{ru}.)

The concordance genus is more elusive and less studied than these two other invariants. 
Gordon \cite[Problem 14]{go2} asked whether $g_4(K) = g_c(K)$ for all knots.  Casson, in
unpublished work,   used the Alexander polynomial to show that the knot $6_2$ satisfies
$g_4(6_2) = 1$ and
$g_c(6_2) = g(6_2) = 2$.  Independently, Nakanishi \cite{na} used a similar argument to
give  examples showing that the gap between $g_4(K)$ and $g_c(K)$ can be arbitrarily large,
for knots with
$g_4$ arbitrarily large.   In Section \ref{cagenus} we  briefly review these results
and give what is essentially Nakanishi's example 
showing that
$g_c(K)$ can be arbitrarily large for knots with $g_4(K) = 1$.  (Obviously, if $g_4(K) = 0
$ then
$g_c(K) = 0$.)  We then show that by using the signature  in conjunction with the
Alexander polynomial we can attain finer results: we   construct   knots $K$ with $g_4(K)
= 2$ and with the same Alexander polynomial as a slice knot, but with  $g_c(K)$ arbitrarily
large.

\sh{Algebraic concordance and higher dimensional knot theory}

Associated to a knot $K$ and choice of Seifert surface, $F$,  there is a Seifert form
$V_K$: this is an integral matrix satisfying det$(V_K - V_K^t) = \pm 1$, where $V_K^t$
denotes the transpose.  There is a Witt group of such Seifert forms, denoted
$\calg_-$, defined by Levine \cite{le1}. Denoting the concordance group of knots by
$\calc_1$, Levine proved that the map 
$K
\to V_K$ induces a homomorphism
$\psi:
\calc_1 
\to
\calg_{-}$.  

Knot invariants that are defined on $\calg_-$ are called {\it algebraic} invariants, and
it is easily shown that the Alexander polynomial and signature based obstructions are
algebraic. A general algebraic invariant of a knot, $g_c^a(K)$, is defined to be 
one half the rank of the minimal dimension representative of $V_K$ in
$\calg_-$.  Everything we have discussed so far generalizes to higher dimensional
concordance, where Levine proved that  the map $\psi$ classifies knot concordance. Hence
we have:

\begin{theorem}  In higher dimensions, $g_c(K) = g_c^a(K)$. \end{theorem}

\noindent (Given a  knot $K$ we can also form the hermitian matrix $(1-z)V_K + (1- 
{z^{-1}})V_K^t$, over the field of fractions of $\qq[z]$, $\qq(z)$.  This induces a well defined
homomorphism
$\psi':\calc_1 \to W(\qq(z))$, where $W(\qq(z))$ is the Witt group of hermitian forms
on vector spaces over the function field $\qq(z)$.  There is an invariant $g_4^a(K)$
given by the minimal rank representative of the class of
$\psi'(K)$.  It can be shown that $g_4(K) \ge g_4^a(K)$ and we conjecture that 
  in higher dimensions this becomes an equality.)

\sh{Algebraically slice knots}

Our deepest and most subtle results concern algebraically slice knots. We begin with a
definition: 
\begin{definition}  The map
$\psi\! : \calc_1 \to \calg_-$ has kernel denoted $\cala$, the {\it
concordance group of algebraically slice knots}.
\end{definition}

Four--dimensional knot concordance is  unique and especially challenging in that, unlike
the  higher dimensional analogs, $\cala$ is nontrivial.  In the smooth setting a number
of techniques based on the work of Donaldson \cite{do} and Witten \cite{wi} (see for
example \cite{km}) have given new insights into the structure of $\cala$.  However, in the
topological locally flat category the  only known obstructions to a knot in $\cala$ being
trivial are Casson--Gordon invariants \cite{cg1, cg2} and  their extensions (for example
\cite{cot}).  In the language of the present paper, the  results of \cite{cg1, cg2} can be
stated as:

\begin{theorem}  There exist knots $K \in \cala$ with $g_4(K) \ge 1.$
\end{theorem}

Gilmer extended this result in \cite{gl1}:

\begin{theorem}  For every $N$ there exist knots $K \in \cala$ with $g_4(K) \ge N$.
\end{theorem}

In the Casson--Gordon examples of twisted doubles of the unknot one has 
that $g_4(K) = g(K) =
1$. In Gilmer's examples  $g_4(K) = g(K) = N$.

Our main result concerning $\cala$ is the following.

\begin{theorem}\label{result} For every $N$ there exists a knot $K \in \cala$ with $g_4(K) =
1$ and
$g_c(K) = g(K) = N$.
\end{theorem}

To conclude this introduction we remark on the inherent challenge of proving Theorem
\ref{result}.  Showing that a given algebraically slice knot is not slice is equivalent to
showing that   it is not concordant to a single knot, the unknot.  In the case, say, of
showing that  a genus 2 algebraically slice knot is not concordant to a knot of genus 1, we have to prove
that it is not concordant to any knot in an infinite family of knots, each of which is algebraically slice and hence about which one
knows very little.  There are of course some constraints on this family of knots based on
their being genus 1, such as the Alexander polynomial, but with the added restriction that the knots are algebraically slice these do not apply to the present problem.  The main remaining tools are Casson--Gordon invariants, however  known genus  constraints based on these~\cite{gi1} already would apply to bound the  4--ball genus as well, so these cannot be directly applicable either.  As we will see, Casson--Gordon invariants still are
sufficient to provide examples, but the proof calls on two steps.  The first is a delicate
analysis of metabolizing subgroups for the linking forms that arise in this problem.  The
second is the construction of knots with Casson--Gordon invariants satisfying   rigid
constraints.

\sh{References and conventions}

We will be working in the smooth category throughout this paper.  All the results carry
over to the topological locally flat category by \cite{fq}.  

 Basic references for knot theory include
\cite{bz, ro}.  The fundamentals of concordance and Levine's work are contained in
\cite{le1, le2}. The principal references for Casson--Gordon invariants are the original
papers,
\cite{cg1, cg2}.

\section{Algebraic bounds on the concordance genus}\label{cagenus} 

In this section we will study bounds on $g_c$ based on the Seifert form of the knot.  All
of these are easily seen to depend only on the algebraic concordance class, and hence are
in fact bounds on $g_c^a(K)$.  Because of this,  none can yield information regarding
$g_c(K)$ for knots $K \in \cala$.

\subsection{Alexander polynomial based bounds on $g_c$} Recall that the Alexander polynomial
of a knot
$K$ is defined to be
$\Delta_K(t) = $ det($V_K - tV_K^t$) where
$V_K$ is an arbitrary Seifert matrix for $K$.  It is well defined up to multiplication by
$\pm t^{n}$ so we will assume that $\Delta_K(t) \in \zz[t]$ and $\Delta_K(0) \ne 0$.  The
degree of such a representative will be called the degree of the Alexander polynomial,
deg($\Delta_K(t))$.  

A simple observation  regarding the Alexander polynomial and concordance is that if a
Seifert form $V$ represents 0 in $\calg_-$ then $\Delta_V(t) = \pm t^n f(t)f(t^{-1})$ for
some polynomial $f$ and integer $n$. (This result was mentioned in \cite{fo} and first
proved in \cite{fm}.)  It follows  that if
$V_1$ and
$V_2$ represent the same class in $\calg$ then $\Delta_{V_1}(t)  \Delta_{V_2}(t) =
\pm t^n  f(t)f(t^{-1})$ for some polynomial $f$.  From this we  have the following basic example of
Casson.

\begin{example}  {\sl The knot $6_2$, illustrated in figure 1, satisfies $g_c(6_2) = 2$ and $g_4(6_2)
= 1$.}  Note first that   $\Delta_{6_2}(t) = t^4 - 3t^3 + 3t^2 - 3t +1$ (\cite{ro}), an
irreducible polynomial. Hence, if $6_2$ were concordant to a knot of genus 
1, we would then have
$\Delta_{6_2}(t)g(t)= \pm t^n f(t)f(t^{-1})$ for some polynomial $f$, integer $n$, and
polynomial $g(t)$ with deg($g(t)) \le 2$. Degree considerations show that this is
impossible. On the other hand, Seifert's algorithm applied to the standard diagram of
$6_2$ yields a Seifert surface of genus 2.

To see that $g_4(6_2) = 1$, observe that the unknotting number of $6_2$ is 1 (change the
middle crossing) and so $6_2$ bounds a surface of genus 1 in the 4--ball.  It follows that
$g_4(K) \le 1$.  On the other hand
$6_2$ is not slice since its Alexander polynomial is irreducible, so it cannot bound a
surface of genus 0.

\end{example}

\cfigure{1.1in}{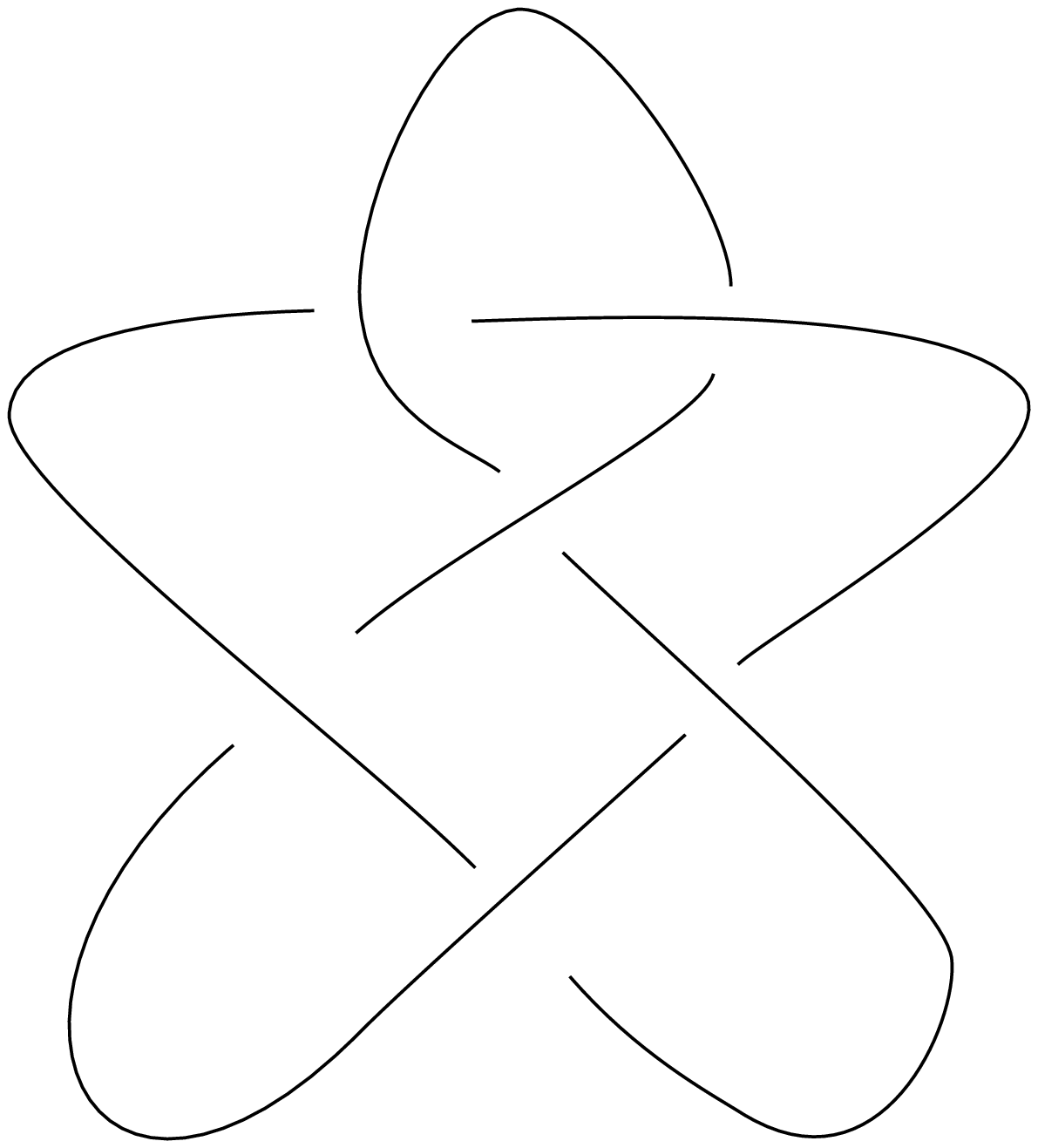}{The knot $6_2$}

Nakanishi \cite{na}, independently of Casson, used the Alexander polynomial in
the same way to develop other examples contrasting $g_c$ and $g_4$.  These techniques are
summarized by the following theorem.

\begin{theorem} \label{thmcgenus}  Suppose  that $  \Delta_K(t)$ has an irreducible
factorization in $\qq[t]$ as   
$$\Delta_K(t) = p_1(t) ^{\epsilon_1} \cdots p_k^{\epsilon_m} q_1(t)^{\delta_1} \cdots
q_j(t)^{\delta_j}$$
 where the $p_i(t)$ are distinct irreducible  polynomials with $p_i(t) = \pm t^{n_i}
p_i(t^{-1})$ for some $n_i$ and  $q_i(t)
\ne \pm t^{n_i}q_i(t^{-1})$ for any $n_i$.  Then $g_c(K) $ is greater than or equal to one
half the sum of  the degrees of the $p_i$ having exponent $\epsilon_i$ odd. 
\end{theorem}

Using this, Nakanishi proved the following.  (In fact, he gives similar examples with other
values of $g_4(K)$.)  We include this argument because a related construction is used in the
next subsection.

 \begin{theorem} For every $N >0$ there exists a knot $K$ with $g_4(K) = 1$ and
$g_c(K) > N$.

\end{theorem}

\begin{proof} According to Kondo and Sakai, \cite{ko, sa}, every Alexander polynomial
occurs as the Alexander polynomial of an unknotting number one knot.  Hence, the proof is
completed by finding irreducible Alexander polynomials of arbitrarily high degree.  Such  
examples include  the cyclotomic polynomials $\phi_{2p}(t)$ with $p$ an odd prime. It is
well known that cyclotomic polynomials are irreducible.  We have that
$$ \phi_{2p}(t) =  { \frac{  (t^{2p} -1)(t-1)} {  (t^2 -1)(t^p -1)  }  } = t^{p-1} - t^{p-2} + t^{p-3} -
\ldots t +1.$$  This is an Alexander polynomial since $\phi_{2p}(t)$ is symmetric and
$\phi_{2p}(1) = 1$. Hence, the unknotting number one knot with this polynomial has $g_4(K)
= 1$ but
$g_c(K)
\ge (p-1)/2$.
 \end{proof}

An examination of the construction used by Sakai  in \cite{sa} shows that the knot used
above also has 
$g(K) = (p-1)/2$.  Briefly, the knot is constructed from the unknot by performing $+1$
surgery in
$S^3$ on an unknotted circle $T$ in the complement of the unknot $U$.  The surgery circle
$T$ meets a disk bounded by $U$ algebraically 0 times but geometrically $(p-1)$ times.
Hence, a genus $(p-1)/2$ surface bounded by $U$ that misses $T$ is easily constructed.

\subsection{Further bounds on $g_c$}

Certainly this inequality of Theorem \ref{thmcgenus} cannot be replaced with an equality.  

\begin{example} {\sl The granny knot (the connected sum of the trefoil with itself) has
concordance  genus 2  and has Alexander polynomial
$(t^2 - t +1)^2$.  The square knot, the connected sum of the trefoil with its mirror image
has the same Alexander polynomial but has concordance genus 0.}  
 To see this, first recall that both these knots have genus 2.  We have that
$g_c(K)
\ge g_4(K)$.  According to Murasugi
\cite{mu}, the classical signature of a knot bounds $g_4$; more precisely,   $g_4(K)
\ge {\frac{1 }{ 2}}
\sigma(K)$.  The signature of the granny knot is 4, and hence we have the desired value of
$g_4$ for the granny knot.  On the other hand, the square knot is of the form $K \# -K$
and hence is slice.
\end{example}

The rest of this subsection   will discuss strengthening Theorem \ref{thmcgenus}.  We
begin by recalling Levine's construction of isometric structures in \cite{le1}.  Every
Seifert form
$V$ is equivalent (in
$\calg_-$) to a nonsingular form of no larger dimension.  Associated  to such a $V$ of
dimension $m$ we have an isometric structure $( \left< , \right> , T)$ on a rational
vector space $X$ of dimension $m$, where $\left< , \right>$ is the quadratic form on
$X$ given by $V + V^t$ and $T$ is the linear transformation of $X$ given by
$-V^{-1}V^t$. The map $V \rightarrow  (V+ V^t, -V^{-1}V^t)$ defines an isomorphism from
the Witt group of rational Seifert forms $\calg^\qq$ to the Witt group of rational
isometric structures, $\calg_\qq$. The Alexander polynomial of $V$ is the characteristic
polynomial of $T$.  (In the Witt group $\calg_\qq$ an isometric structure is Witt trivial,
by definition, if the inner product $  \left< , \right>  $ vanishes on a half--dimensional
$T$--invariant subspace of $X$.)

As a $\qq[t,t^{-1}]$ module $X$ splits as a direct sum $\oplus X_{p(t)}$ over all
irreducible polynomials $p(t)$, where $X_{p(t)}$ is annihilated by some power of $p(t)$. 
 According to Levine, any isometric structure is equivalent to one with the
$X_{p(t)}$ trivial if $p(t) \ne \pm t^n p(t^{-1})$ for some $n$. Furthermore,
\cite[Lemma 12]{le2}, each remaining
$X_{p(t)}$ can be reduced to a  Witt equivalent  form annihilated by $p(t)$: $$X_{p(t)} = 
\left({   \frac{\qq[t,t^{-1}]
}{ <p(t)>  }    }\right)^k $$ for some $k$. 

Write $X $ as $\oplus_{i=1 \ldots s} X_{p_i}$ where the $X_{p_i}$ are all of the given
form.    Now, suppose that
$p_i(t)$ has as a root
$e^{i
\theta}$ for some real
$\theta$.  The Milnor 
$\theta$--signature of $V$, $\sigma_\theta(V)$, (see \cite{mi}) is defined to be the
signature of the quadratic form
$\left< ,
\right>$ restricted to the ({\it real\/}) summand of $X_{p_i(t)}$ associated to $p_\theta(t)
= t^2 - 2
\cos(\theta) t +1$.  From this analysis the  next theorem follows immediately.

\begin{theorem} Suppose $\Delta_V(t)$ has distinct symmetric irreducible factors
$p_i(t)$ and
$p_i(e^{i
\theta_i }) = 0$.  If $\sigma_{\theta_i}(V) = 2k_i$ then $g^a_c(V) \ge {\frac{1 }{ 2}}\sum_{i}
|k_i|({\rm deg} (p_i))  $.

\end{theorem}

\noindent Notice that there can be distinct values of $\theta_i$ for which $p_i(e^{i
\theta_i }) = 0$.  

In general, the computation of the
Milnor
$\theta$--signatures can be nontrivial.  The following examples illustrate how the
signature  used in conjunction with the Alexander polynomial yields much stronger results
than can be obtained using either one alone.

\begin{example} For a given   prime $p = 3 \mod 4$, consider an unknotting number 1 knot
$K$ with
$\Delta_K(t) =  \phi_{2p}(t) $.  According to Murasugi \cite{mu}, if $|\Delta_K(-1)| = 3
\mod 4$ then 
$\sigma(K) = 2
\mod 4$, where $\sigma(K)$ is the classical knot signature, the signature of $V_K + V_K^t$. It
is easily shown that
$\phi_{2p}(-1) = p$, so for our
$K$ we have
$|\sigma(K)| = 2 \mod 4$.  However, since $g_4(K) = 1$, $|\sigma(K)| \le 2$.  After changing
orientation if need be, we have that $\sigma(K) = 2$.  By \cite{ma} $\sigma(K)$ is given as a
sum of Milnor signatures, so it follows that for some
$\theta$,
$\sigma_\theta(K) = 2$.  Now, let
$J = K
\# K$.  Since Milnor signatures are additive under connected sum, $\sigma_\theta(J) = 4$.
We also have
$\Delta_J(t) = (\phi_{2p}(t))^2$, which is  of degree
$2p-2$.  Hence by the previous theorem, $g_c(J) \ge \mbox{deg}(\phi_{2p}(t)) = p-1$.  No
bound on
$g_c$ can  be obtained using Theorem \ref{thmcgenus} since the polynomial is a square. 
Since $J$ is unknotting number 2, we have $c_4(J) \le 2$ but the signature implies that
$g_4(J) = 2$. 

\end{example}

\section{Casson--Gordon invariants}

\subsection{Basic theorems} We will be working with a fixed prime number $p$ throughout
the following discussion.
 
For a knot $K$ let $M(K)$ denote the 2--fold branched cover of $S^3$ branched over
$K$.  Let $H_K$ denote the $p$--primary summand of $H_1(M(K),\zz)$.  More formally, we have
$H_K = H_1(M(K), \zz_{(p)})$, where $\zz_{(p)}$ represents the integers localized at
$p$; in other words, $\zz_{(p)} = \{ {\frac{m}{ n}} \in \qq | \gcd (p,n) = 1\}$.

There is a nonsingular symmetric linking form   $\beta\! : H_K \times H_K \to
\qq/\zz$. If
$K$ is algebraically slice there is a subgroup $M \subset H_K$ satisfying $M = M^\perp$
with respect to the linking form.  Since the linking form is nonsingular, this easily
implies that $|M|^2 = |H_K|$.  Such an $M$ is called a {\it metabolizer} for $H_K$.

 Let
$\chi\!: H_K \to \zz_{p^k}$ be a homomorphism.  The Casson--Gordon invariant 
$\sigma(K,\chi)$ is a rational invariant of the pair
$(K,\chi)$.  (See \cite{cg1}, where this invariant is denoted $\sigma_1
\tau(K,\chi)$  and $\sigma$ is used for a closely related invariant.)  The main result in
[CG1] concerning Casson--Gordon invariants and slice knots that we will be using is the
following.

\begin{theorem} \label{CG} If $K$ is slice then  there is a metabolizer $M \subset H_K$ such
that 
$\sigma(K,\chi) = 0$ for all 
$\chi\! : H_K \to \zz_{p^k}$ vanishing on $M$.\end{theorem}

 We will be using  Gilmer's additivity theorem  \cite{gi2}, a vanishing result proved
by Litherland \cite[Corollary B2]{lit2}, and a simple fact that follows immediately from
the definition of the Casson--Gordon invariant.

\begin{theorem}  \label{Gi} If $\chi_1$ and $\chi_2$ are defined on $M_{K_1}$ and  
$M_{K_2}$, respectively, then  
$\sigma( K_1\ \#\ K_2, \chi_1\ \oplus \
\chi_2) =
\sigma(K_1,\chi_1) + \sigma(K_2,\chi_2)$. \end{theorem}

\begin{theorem} If $\chi$ is the trivial character, then $\sigma(K,\chi) = 0$.
\end{theorem}
\begin{theorem} For every character  $\chi$, 
$\sigma(K,\chi) = \sigma(K,-\chi)$.
\end{theorem}

\subsection{Identifying characters with metabolizing elements}  We will be considering
characters that vanish on a given metabolizer $M$.  Note that the character given by
linking with an element  $m \in M$ is such a character and that any character 
$\chi\! : H_K \to \zz_{p^k}$  vanishing on $M \subset H_K$ is of the form
$\chi(x) = \beta(x,m)$ for some
$m \in M$.   We will denote this character by $\chi_m$.

\subsection{Companionship results} \label{companions} Our construction of 
examples of
algebraically slice knots will begin with a knot $K$ with a null homologous link of $k$
components in the complement of $K$, $L = \{L_1, \ldots ,L_k \}$.  $L$ will be an unlink,
though it will link $K$ nontrivially.  A new knot,
$K^*$, will be formed by removing from $S^3$ a neighborhood of $L$ and replacing each
  component  with the complement of a knot, $J_i$.  This
can be done in such a way that the resulting manifold is again $S^3$.  (The attaching map
should identify the meridian of
$L_i$ with the longitude of $J_i$ and vice versa.) The image of $K$ in this new copy of
$S^3$ is the knot we will denote $K^*$.  

Let $\tilde{L}_i$ denote a lift of $L_i$ to the 2--fold branched cover, $M(K)$.  There is
a natural identification of $H_1(M(K),\zz)$ and $H_1(M(K^*),\zz)$.  Suppose that $\chi\! :
H_K
\to \zz_{p^j}$ and that $\chi(\tilde{L}_i) = a_i$.  We have the following theorem relating
the associated Casson--Gordon invariants of $K$ and $K^*$.  A proof is basically contained
in \cite{gl1}.  The result is implicit in \cite{lit2} and \cite{gi2}.  In the formula,
$\sigma_{a_i/p^j}(J_i)$ denotes the classical Tristram--Levine signature \cite{tr} of $J_i$. 
This signature is defined to be the signature of the hermitian form 
$$(1 - e^{ {\frac{a_i}{ p^j }} 2 \pi i}) V_{J_i} + (1 - e^{ -{\frac{a_i }{ p^j}} 2 \pi i}) V^t_{J_i}.$$

\begin{theorem}\label{companiont}  In the setting just described, 
$$\sigma(K^*,\chi) -
\sigma(K,\chi) = 2\sum_{i=1}^k \sigma_{a_i/p^j}(J_i).$$

\end{theorem}

\section{Properties of metabolizers}

In the next section we will construct an algebraically slice  knot $K$ with $g(K) = N$ and
$H_K
\cong (\zz_3)^{2N}$.  We will show that it is  not concordant to a knot $J$ with
$g(J) < N$ by proving that if $K$ is concordant to $J$ then rank($H_J) \ge 2N$.  The
following is our main result relating Casson--Gordon invariants and genus.  With the exception of one example our applications all occur in the case of $p=3$.

\begin{theorem} \label{met4}  If $K$ is an algebraically slice knot with $H_K
\cong (\zz_p)^{2g}$ and
$K$ is  concordant to a knot $J$  of genus $g' < g$ then there is a metabolizer $M_K
\subset H_K$ and a nontrivial subgroup $M_0 \subset M_K$ such that for $\chi_m$ with $m
\in M_K$, 
$\sigma(K,\chi_m)$ depends only on the class of $ m$ in the quotient
$  M_K / M_0$.  That is, if $m_1 \in M_K$ and $m_2 \in  M_K$ with $m_1 - m_2 \in M_0$, then
$\sigma(K ,
\chi_{m_1}) =  \sigma(K  ,
\chi_{m_2})$.
\end{theorem}

\subsection{Metabolizers} \label{metabolizers}
\begin{theorem}\label{met1thm}  If $K$ is an algebraically slice knot of genus $g$, the
linking form on $H_K$ has a metabolizer generated by $g$ elements.

\end{theorem}

\begin{proof}  Because $K$ is algebraically slice, with respect to some generating set its
Seifert matrix is of the form
$$\left( \begin{matrix} 0 & A \\
               B & C
\end{matrix} \right) $$ for some $g \times g$ matrices $A$, $B$, and $C$.  Hence,
$H_1(M(K),\zz)$ has homology presented by $V_{K } + V_{K }^t  $, which is of the form
$$P = \left( \begin{matrix} 0 & D \\
               D^t & E
\end{matrix} \right) $$ for other matrices, $D$ and $E$, where $D$ has nonzero
determinant. The order of
$H_1(M(K),\zz)$ is det$(D)^2$.  

This presentation matrix corresponds to a generating set $\{x_i, y_i\}_{i = 1, \ldots ,
N}$.  We claim that the set $\{ y_i\}$ generates a metabolizer.  First, to see that it is
self--annihilating with respect to the linking form, we recall that with respect to the
same generating set the linking form is given by the matrix
$$P^{-1} = \left( \begin{matrix} -(D^{-1})^tE D^{-1} &  (D^{-1})^t\\
               D^{-1} & 0
\end{matrix} \right). $$ That this is the correct inverse can  be checked by direct
multiplication.  The lower right hand block of zeroes implies the vanishing of the linking
form on $< \{y_i\}>$.

We next want to see that  $\{y_i\}$ generate a subgroup of order det$(D)$.  Clearly the
$ y_i $ satisfy the relations given by the matrix $D$.  What is not immediately clear is
that the relations given by $D$ generate all the relations that the $\{y_i\}$  satisfy. 
To see this, note that any relations satisfied by the $\{y_i\}$ are given as a linear
combination of the rows of $P$.  But since the block $D$ has nonzero determinant, any such
combination will involve the $\{x_i\}$ unless all the coefficients corresponding to the
last $g$ rows of $P$ vanish.  This implies that the relation comes entirely from the
matrix $D$.
\end{proof}

\noindent{\bf Notation}\qua Suppose that the algebraically slice knot $K$ is concordant to a
knot $J$.  Let
$M_\#$ be a metabolizer for   $H_{K \# -J}$.  Let $M_J$ be a metabolizer for $H_J$.  Let
$M_K = \{ m \in H_K\ | \  (m,m') \in M_\# \mbox{ for some } m' \in M_J\}$.  For each
element $m' \in M_J$, set $M_{m'} = \{ m \in H_K\ |
\  (m,m')
\in M_\# \}$.  In particular, $M_{0} =   \{ m \in H_K\ | \  (m,0)
\in M_\#\}$ and $  M_K = \cup_{m' \in M_J} M_{m'}$. Finally, let $M_{J,0} =
\{ m'
\in H_J\ |
\  (0,m')
\in M_\# \}$.

\begin{theorem}\label{met2thm}  With the above notation, $M_K$ is a metabolizer for
$H_K$. 
\end{theorem}

\begin{proof} A proof of the corresponding theorem for bilinear forms on vector spaces
appears in
\cite{ke}.  A parallel proof for finite groups and linking forms can be constructed in a
relatively straightforward manner.  One such proof appears in \cite{klv}.  Since all
metabolizers split over the $p$--primary summands, the results follow for these summands.
\end{proof}

The set of elements $M_{0} $  is surely nonempty: it contains 0.  It is also easily seen
to be a subgroup.  

\begin{lemma} If $M_{m'}$ is nonempty then it is a coset of $M_0$ in $M_K$.\end{lemma}

\begin{proof} The proof is straightforward.  If $x, y \in M_{m'}$ then
$(x,m') \in M_{\#}$ and 
$(y,m') \in M_{\#}$.  Hence, $(x-y,0) \in M_\#$, so $x-y \in M_0$. Similarly, if  $x
\in M_{m' }$ and $y \in M_0$, then $(x,m' ) \in M_\#$ and 
$(y,0) \in M_{\#}$, so $(x+y,m') \in M_{\#}$ and 
$ x+y  \in M_{m'}$.
\end{proof}

\begin{lemma}\label{met3thm} The map $M_{m'} \rightarrow m'$ induces an injective
homomorphism  of $M_K/M_0$ to $ M_J / M_{J,0}$.
\end{lemma}
\begin{proof} It must be checked   that this map is well-defined.  Suppose first that
$M_{m'} = M_{m''}$.  Then for any $m \in M_{m'} = M_{m''}$, $(m,m') \in M_\#$ and $(m,m'') \in
M_\#$.  Taking differences, we have that $(0,m' - m'' )\in M_\#$, implying that $m' - m'' \in
M_{J,0}$ as desired.

That this map is a homomorphism is trivially checked.

To check injectivity, we need to show that for all $m' \in M_{J,0}$, $M_{m'} = M_0$.  But
$0 \in M_{m'}$ since $(0,m') \in M_\#$ by the definition of $M_{J,0}$.  Since $0
\in M_{m'}$, $ M_{m'}$ is the identity coset, as needed.
\end{proof}

\begin{theorem} Let $K$ be an algebraically slice knot with $H_K \cong (\zz_p)^{2g}$ and
suppose that $K$ is concordant to a knot $J$ with $g(J) < g$. Then for some metabolizer
$M_\#$ for $H_{K \# - J}$ and for any metabolizer $M_J$ for $H_J$, the subgroup 
$M_0 \subset H_K$ is nontrivial.
\end{theorem}

\begin{proof} If $M_0$ is trivial we would have, by Lemma \ref{met3thm}, an injection of
$(\zz_p)^{g}$ into $M_J / M_{J,0}$.  But by Theorem \ref{met1thm} the metabolizer
$M_J$ can be chosen so that it has rank less than $g$.  It follows that a quotient will
also have rank less than $g$. Hence, it cannot contain a subgroup of rank $g$.
\end{proof}

We now prove Theorem \ref{met4}:

\medskip

\noindent{\bf Theorem \ref{met4}}\qua {\sl  If $K$ is an algebraically slice knot with $H_K
\cong (\zz_p)^{2g}$ and
$K$ is  concordant to a knot $J$  of genus $g' < g$ then there is a metabolizer $M_K
\subset H_K$ and a nontrivial subgroup $M_0 \subset M_K$ such that for $\chi_m$ with $m
\in M_K$, 
$\sigma(K,\chi_m)$ depends only on the class of $ m$ in the quotient
$  M_K / M_0$.  That is, if $m_1 \in M_K$ and $m_2 \in  M_K$ with $m_1 - m_2 \in M_0$, then
$\sigma(K ,
\chi_{m_1}) =  \sigma(K  ,
\chi_{m_2})$.}

\begin{proof} Since $K \# -J$ is slice, we let $M_\#$ be the metabolizer given by Theorem
\ref{CG}. We also have that 
$-J$ is algebraically slice, so we let $M_J$ be an arbitrary metabolizer for $H_J$ with
rank$(M_J) < g$ and we let $M_K \subset H_K$ be the metabolizer constructed above.  We also
let
$M_0$ be the nontrivial subgroup of $M_K$ described above.

Let $\chi_{m_1}$ and $\chi_{m_2}$ be characters on $H_K$ vanishing on $M_K$.   We are
assuming  further that $m_1$ and
$m_2$ are in the same coset of $M_0$: $m_1$ and $m_2$ are both in $M_{m'}$ for some $m'
\in M_J$.  We want to show that
$\sigma(K,\chi_{m_1}) =
\sigma(K,\chi_{m_2})$. 

Since $m_i \in M_{m'}$,  we have  that $(m_1,m')
\in  M_\# $ and $(m_2,m')\in  M_\#$.  Hence, by Theorem \ref{CG}, $$\sigma(K \# -J,
\chi_{m_1} \oplus \chi_{m'} )= 0 = \sigma(K \# -J,
\chi_{m_2} \oplus \chi_{m'} )$$ The result now follows immediately from the additivity of
Casson--Gordon invariants.
\end{proof}

\section{Construction of examples}

\subsection{Description of the starting knot, $K$}  We will build a knot $K^*$ with the
desired properties regarding $g_c$.      The construction begins with a knot $K$ which is
then modified to build $K^*$.  In this subsection we describe $K$ and its properties.

Figure 2 illustrates a knot $K$ and a link $L$ in its
complement.   The figure is drawn for the case $N = 3$.  The correct generalization for
higher $N$ is clear.  Ignore
$L$ for now.     The knot
$K$ bounds an obvious Seifert surface
$F$ of genus N.  The Seifert form of $K$ is $$N\left( \begin{matrix} 0 & 1 \\
               2 & 0
\end{matrix} \right). $$ The homology of $F$ is generated by the
symplectic basis $\{x_i, y_i\}_{i = 1, \ldots, N}$. Here each  $x_i$ is represented by
the simple closed curve formed as the union of an embedded arc  going over the left band of
$i^{th}$ pair of bands and an embedded arc in the complement of the set of bands.  The
$y_i$ have similar representations, using the right side band of each pair. 

The knot
$K$ is assured to be slice by arranging that the link formed by any collection
$\{z_i\}_{i = 1, \ldots, N}$, where each $ z_i$ is either $x_i $ or $y_i$, forms an
unlink.  (We are not distinguishing here between the class  $x_i$ and the embedded curve
representing the class; similarly for
$y_i$.) 

\cfigure{5in}{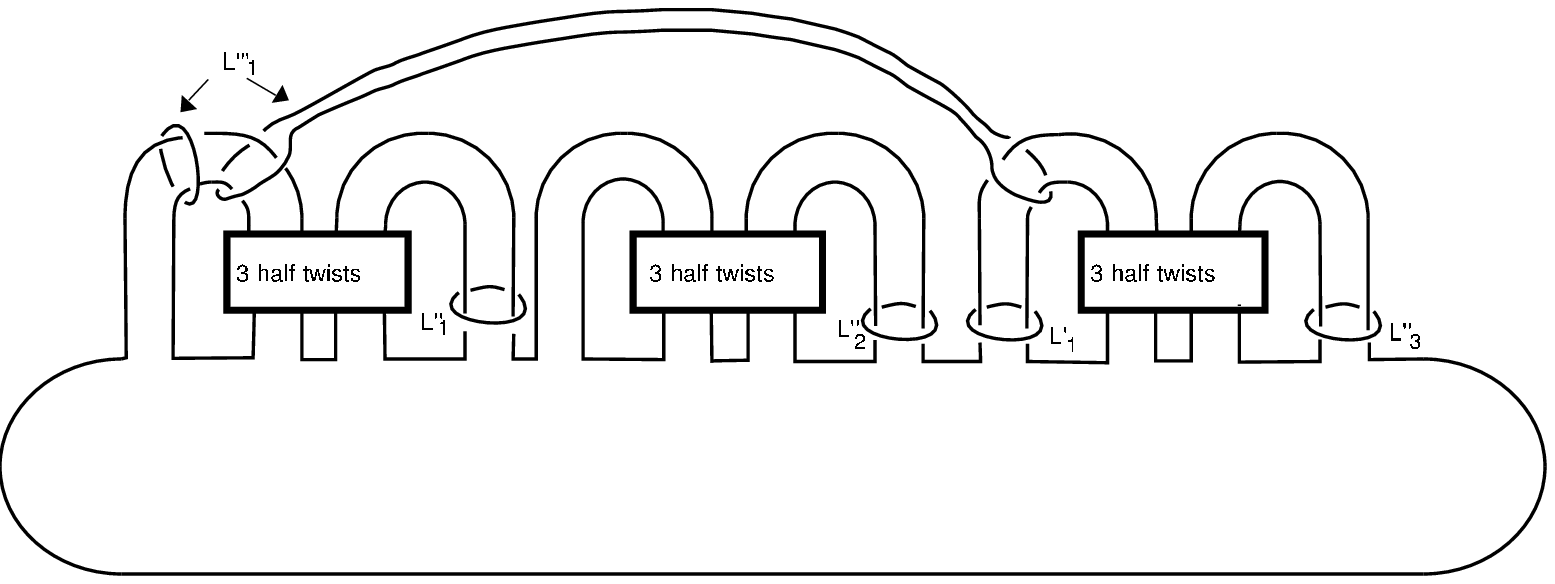}{The basic knot}

 The homology of the complement of $F$ is generated by trivial linking curves to the
bands, say
$\{a_i, b_i\}_{i = 1 , \ldots, N}$. The 2--fold cover of $S^3$ branched over $K$, $M(K)$,
satisfies
$H_1(M(K),\zz) \cong (\zz_3)^{2N}$.  Picking arbitrary lifts of the $\{a_i, b_i\}_{i = 1 ,
\ldots, N}$ gives a set of curves in $M(K)$,  $\{\tilde{a}_i,
\tilde{b}_i\}_{i = 1 ,
\ldots, N}$, generating $H_1(M(K),\zz)$.  This follows from
standard knot theory constructions \cite{ro}, but perhaps is most evident using the surgery
description of $M(K)$ given by Akbulut and Kirby \cite{ak}.  It also follows easily from
this description of $M(K)$ that the linking form with respect to $\{\tilde{a}_i,
\tilde{b}_i \}_{i = 1, \ldots, N}$ is given by
$$N\left( \begin{matrix} 0 & {\frac{1}{  3}} \\
              {\frac{1}{  3}} & 0
\end{matrix} \right). $$
 Here there is a slight issue of signs, but the signs as given in this linking matrix can
be achieved by choosing the appropriate lifts,  or simply by orienting the lifts properly.

\subsection{Construction of the link $L$}

The desired knot $K^*$ is constructed from $K$ by removing the components of a link $L$ in
the complement of $F$ and replacing them with complements of knots $J_i$.  In this
subsection we describe $L$.

The link $L$ consists of three sublinks: $L = \{L', L'', L'''\}$.  Here is how the various
components of $L$ are chosen:

\begin{itemize}
\item $L'$ has only one component: $L' = \{L'_1\}$.  Here $L_1$ is chosen to be a trivial
knot representing
$a_N$, the linking circle to the band with core $x_N$.

\item $L''  = \{L''_i\}_{i =1, \ldots , N }$.  We choose $L''_i$  to be a trivial knot
representing
$b_i$, the linking circle to the band with core $x_i$.

\item $L''' = \{L'''_i\}$ consists of a set of  2--component sublinks.  For each ordered
pair, $( a_i , b_j)_{ i = 1, \ldots , (N-1) , j = 1, \ldots , N}$ we have a two component
link $L'''_i$:  one component is  a trivial knot representing $a_i$ as a small linking
circle to $x_i$; the other component is the band connected sum of a curve parallel to that
one with a trivial knot representing $b_j$ as a small linking circle to
$b_j$. Similarly, 2--component links are formed for the pairs $(a_i, a_N)$, $ i < N$.  The
set
$L'''$ has 
$N^2 - 1$ elements.

\end{itemize}

 In the figure we have indicated all the components of $L'$ and $L''$.  The only sublink
of $L'''$ that is illustrated is the one corresponding to the ordered pair $(a_1, a_3)$.

This collection is chosen so that the following theorem holds.
\begin{theorem}  

{\bf A}\qua In $S^3 - \{x_i\}_{i = 1, \ldots, (N-1)}$ the components of $L'$ and $L''$ form
an unlink, split from the link  $L''' \cup  \{x_i\}_{ i = 1,
\ldots , (N-1)}$.

{\bf B}\qua  The link  $L''' \cup  \{x_i\}_{ i = 1,
\ldots , (N-1)}$ is the union of an unlink,   $\{x_i\}_{ i = 1,
\ldots , (N-1)}$    with  parallel pairs of meridians to the $x_i$, one pair for each
sublink $L'''_i$ of $L'''$.

 \end{theorem}

\subsection{Constructing $K^*$ and its properties}

 We will be selecting sets of knots $\{J'_i\}$, 
$\{J'_i\}$, and $\{J'''_i\}$.  There is only one knot in the set $\{J'_i\}$; it corresponds
to the knot
$L'_1$. There are $N$ knots in the set $\{J''_i\}$, with one knot $J''_i$ for each
$L''_i$. Finally, there is one knot
$J'''_i$ for each  2--component sublink $L'''_i$ of $L'''$. The necessary properties  of
all these knots will be developed later.   To construct
$K^*$ we follow the companionship construction described in Section  
\ref{companions}:
remove tubular neighborhoods of each
$L'_i$ and
$L''_i$ and replace them with the complement of the corresponding
$J'_i$ or $J''_i$.   Neighborhoods of the two components of
$L'''_i$ are replaced with the complements of the corresponding  $J'''_i $ and its mirror
image, $-J'''_i$.

Since $K^*$ is formed by removing copies of  $S^1 \times B^2$ from $S^3$ and replacing
them with three manifolds with the same homology, the Seifert form of $K^*$ is the same as
that of $K$.   Hence, as for the knot $K$, $H_1(M(K^*),\zz) \cong (\zz_3)^{2N}$ is
presented by 
$$N\left( \begin{matrix} 0 & 3 \\
               3 & 0
\end{matrix} \right).$$ Similarly, the  linking form with respect to the same basis is
presented by the inverse of this matrix,

$$N\left( \begin{matrix} 0 & {\frac{1}{  3}} \\
              {\frac{1}{  3}} & 0
\end{matrix} \right). $$ For framed link diagrams of these spaces, see \cite{ak}.

\subsection{The concordance genus of $K^*$}
 Before proving that the concordance genus of $K^*$ is $N$, we observe the following.

\begin{theorem}  The knot $K$ just constructed has $g(K) = N$ and $g_4(K) = 1$.
\end{theorem}

\begin{proof}   It is clear that $g(K) \le N$.  However,
since the rank of
$H_1(M(K),\zz)$ is
$2N$,
$g(K) \ge N$.

We must now show that
$g_4(K) = 1$.  This is based on the observation that the curves $\{x_i\}_{i = 1 \ldots, {N-1}}$
  form a strongly slice link:  That is, they bound disjoint disks in $B^4$.  To see
this, note that by replacing the components of the
$L'''_i$ with copies of the complements of $J'''_i$ and $-J'''_i$, we have arranged that the
$x_i$  
  have become the connected sums of pairs of the form $J'''_i \# -J'''_i$, and such a
connected sum is a slice knot.  

To build a genus 1 surface in the 4--ball bounded by $K$, simply surger  the Seifert
surface using these slicing disks.
\end{proof}

Since $H_1(M(K^*),\zz) \cong (\zz_3)^{2N}$, we have, in the notation
of Section \ref{metabolizers}, $M_{K^*} \cong (\zz_3)^{ N}$.  We will
be assuming that $K ^*$ is concordant to a knot of lower genus, so we
also have $M_0$ is a nontrivial subgroup of $M_{K^*}$, by Theorem
\ref{met4}.  The proof that $g_c(K^*) = N$ will consist of showing
that the knots $J'_i$, $J''_i$ and $J'''_i$ can be chosen so that the
Casson--Gordon invariants cannot be constant on the cosets of $M_0$.

The following result is a consequence of  Theorem \ref{companiont}.  Notice that there is only one term in the first sum since the link $L'$ has just one component, which links $a_N$.
\begin{theorem} \label{values} $\sigma(K^*,\chi) =  \sigma(K,\chi)\ +\ 2 \sum_i
c_i\sigma_{1/3}(J'_i)\ +\ 2
\sum_i d_i\sigma_{1/3}(J''_i)\  +\  2\sum_i (e_i - e'_i)\sigma_{1/3}(J'''_i)$, where:  

\begin{enumerate}
\item $c_i$ is 0 or 1 depending on whether $\chi(\tilde{a}_N)$ is 0 or not.

\item $d_i$ is 0 or 1 depending on whether $\chi(\tilde{b}_i)$ is 0 or not.

\item The values of the $e_i$ and $e'_i$ are determined as follows.  The element 
$L'''_i  \in  L'''$  corresponds to the class of the form 
$ {a}_k + x
\in H_1(S^3 - F,
\zz)$, where $1 \le k \le N-1$ and either $x= b_l , 1 \le l \le N$, or $x = a_N$.
With this, $e_i$ is 0 or 1 depending on whether $\chi(\tilde{a}_k)$ is 0 or not; $e'_i$ is
0 or 1 depending on whether $\chi(\widetilde{a_k +x})$ is 0 or not. 

\end{enumerate}
\end{theorem}

\noindent{\bf Notation}\qua  If the character $\chi$ is given by linking with an element $m \in
H_{K^*}$, (that is, if $\chi = \chi_m$), then the coefficients $c_i$, $d_i$, and $e_i -
e'_i$ are functions of $m$.  We denote these functions by $C_i$, $D_i$, and $E_i = e_i -
e'_i$.

\begin{theorem} The knots $J'_i$, $J''_i$, and $J'''_i$ can be chosen so that  
$\sigma(K^*, \chi_{m_1}) = \sigma(K^*, \chi_{m_2})$   if and only if the functions $C_i$,
$D_i$, and $E_i$   all agree on $m_1$ and $m_2$.  
\end{theorem}

\begin{proof} The difference $\sigma(K^*, \chi_{m_1}) -  \sigma(K^*, \chi_{m_2})$ is given
by:
\begin{eqnarray*} 
 &\sigma(K , \chi_{m_1}) -  \sigma(K , \chi_{m_2}) &  
  + 2\   \sum   ( C_i(m_1) -  C_i(m_2)) \sigma_{1/3}(J'_i) \\ &&+ 2\   \sum  ( D_i(m_1) - 
D_i(m_2)))  \sigma_{1/3}(J''_i) \\  &&+ 2\  \sum   ( E_i(m_1) -  E_i(m_2))
\sigma_{1/3}(J'''_i) 
\end{eqnarray*} The set of values of $\{|\sigma(K , \chi_{x}) -  \sigma(K , \chi_{y})|\}_{
x, y
\in H_K}$ is  a finite set, so is bounded above by a constant $B$.  Pick $J'_1$ so that
$\sigma_{1/3}(J'_1) > 2B$.  Pick $J''_1$ so that $\sigma_{1/3}(J''_i) >
2\sigma_{1/3}(J'_1)$.  Finally pick each following $J''_i$ and $J'''_i$ so that at each
step the $1/3$ signature has at least doubled over the previous choice.

With this choice of knots the claim follows quickly from an elementary arithmetic
argument. 
\end{proof}

We will now assume that the knot $K^*$ has been constructed using such   collections,
$\{J'_i\}$, 
$\{J'_i\}$, and $\{J'''_i\}$ as given in the previous theorem.  

\begin{lemma} The subgroup $M_0$ must be contained in the subgroup generated by  
$\{\tilde{b}_i\}_{i = 1, \ldots , N-1}$.

\end{lemma}

\begin{proof} Consider a $\chi_m$ with $m \in M_0$.  Write $m$ as a linear combination of
the $\tilde{a}_i$ and $\tilde{b}_i$.  If  $\tilde{b}_N$ or some  $\tilde{a}_i $ has a
nonzero coefficient, then $\chi_m$ will link nontrivially with $\tilde{a}_N$ or some
$\tilde{b}_i$.  In this case, either $C_1(m)$ or some $D_i(m)$ will be nontrivial. (Recall
that the $\tilde{a}_i$ and $\tilde{b}_i$ are duals with respect to the linking form.)
\end{proof}

The proof of Theorem \ref{result}  concludes with the following.

\begin{theorem}  It is not possible for  $\sigma(K,\chi_m)$ to be constant on each coset of
$M_0$. 

\end{theorem}

\begin{proof} To prove this, we have seen that we just need to show that one of the
coefficients, either a
$C_i$, $D_i$ or   $E_i$, is nonconstant on some coset.  In the previous proof we used
the $C_i$ and $D_i$. We now  focus on the
$E_i$.

Using the previous lemma, without loss of generality we can assume that $M_0$ contains an
element $m_0 = \tilde{b}_1 + \sum_{i= 2, \ldots, N-1}  r_i \tilde{b}_i$ for some set of coefficients $r_i$.   

The metabolizer $M_K$ is of order $3^N$, so it must contain  an element not in the span of
$\{\tilde{b}_i\}_{i = 1 \ldots , N-1}$. Adding a multiple of $m_0$ if need be, we can hence
assume that $M_K$ contains an element $m = \tilde{b}_1 + \sum_{i = 2, \ldots , N-1}
\beta_i \tilde{b}_i + \sum_{i = 1 , \ldots , N} \gamma_i \tilde{a}_i + \delta_N
\tilde{b}_N$, with some $\gamma_i$ or $\delta_N$ nonzero.  In fact, by changing sign, and
adding a multiple of $m_0$, we can assume that one of the nonzero coefficients is 1.

We can now select an element from the set $\{\tilde{b}_1, \ldots , \tilde{b}_n,
\tilde{a}_N\}$ on which $\chi_m$ evaluates to be 1.  Denote that element $\tilde{b}$.

We consider the  $L'''_i$ representing the pair
   $\tilde{a}_1$ and $\tilde{a}_1 + \tilde{b}$. In this case we have the following:

\begin{itemize}
\item $\chi_m(\tilde{a}_1) = 1$
\item $\chi_m(\tilde{a}_1 + \tilde{b}) = 2$
\item $\chi_{m-m_0}(\tilde{a}_1) = 0$
\item $\chi_{m-m_0}(\tilde{a}_1 + \tilde{b}) = 1$

\end{itemize}

 Using Theorem \ref{values} we have, for the corresponding $e_i$ and $e'_i$ that:

\begin{itemize}
\item $e_i(\chi_m)(\tilde{a}_1) = 1$
\item $e'_i(\chi_m)(\tilde{a}_1 + \tilde{b}) = 1$
\item $e_i(\chi_{m-m_0})(\tilde{a}_1) = 0$
\item $e'_i(\chi_{m-m_0})(\tilde{a}_1 + \tilde{b}) = 1$
\end{itemize}

Finally, from the definition of $E_i = e_i - e'_i$ we have that $E_i(\chi_{m}) = 0$ and
$E_i(\chi_{m-m_0}) = -1$.  Hence, the Casson--Gordon invariants
cannot be constant on the coset and the proof is complete.
\end{proof}

\section{Enumeration}

We conclude by tabulating the concordance genus for all prime knots with 10 
crossings.  We are also able to compute the concordance genus of all prime knots with fewer than 10 crossings  with two exceptions, the knots $8_{18}$ and $9_{40}$. This is the first such listing.  

In doing this enumeration we have used the knot tables contained in \cite{bz} and
especially the listings of various genera  in
\cite{kaw}.  Results on concordances between low crossing number knots were first compiled
by Conway \cite{co}, and we have also used corrections and explications for Conway's
results taken from \cite{ta}.  In addition, Conway apparently failed to identify three such
concordances---$10_{103}$   and $10_{106}$ are both concordant to the trefoil, $10 _{67}$
is concordant to the knot $5_2$---and those concordances are described below.  We also use the
fact that for all knots $K$ with 10 or fewer crossings, the genus of $K$ is given by half the
degree of the Alexander polynomial.

\medskip
\noindent{\bf Summary}\qua  In brief, there are 250 prime knots with crossing number less than
or equal to 10.  Of these, 21 are slice, and hence $g_c = 0$.  For 210 of them the Alexander
polynomial obstruction yields a bound equal to the genus, and hence for these $g_c = g$. 
There are 17 of the remaining knots which are concordant to lower genus knots for which $g_c$
is known. Finally, there are two knots $8_{18}$ and $9_{40}$, for which $g_3 = 3$ but for which
we have not been able to show that $g_c = 3$.  In the smooth category both of these have $g_4 = 2$ (see for instance the table in \cite{sh}) and in the topological category  $g_4$ seems to be unknown for both, being either 1 or 2.

\subsection{Slice knots}

Among prime knots of 10 or few crossings there are 21  slice knots.  These are
$6_1$, $8_8$, $8_9$, $8_{20}$, $9_{27}$, $9_{41}$, $9_{46}$, $10_{3}$, 
 $10_{22}$,  $10_{35}$, $10_{42}$,  $10_{48}$,  $10_{75}$,  $10_{87}$,  $10_{99}$, 
$10_{123}$,  $10_{129}$,  $10_{137}$,  $10_{140}$,  $10_{153}$, and  $10_{155}$.

\subsection{Examples for which the polynomial condition suffices}

For 210 of the 250 prime knots of 10 or fewer crossings, Theorem
 \ref{thmcgenus} gives a bound for $g_c(K)$ which is equal to $g(K)$.  Hence for all these
knots we have 
$g_c(K) =
\mbox{deg}(\Delta_K(t)) / 2$. This leaves 19 knots.  These are:
\begin{itemize}
\item  $8_{10} , 8_{11} ,8_{18}$

\item $ 9_{24} ,  9_{37},9_{40}$

\item $10_{21}  , 10_{40}  ,  10_{59} ,   10_{62}, 10_{65} , 10_{67},10_{74}, 10_{77},
10_{98}   ,   10_{103} ,
10_{106} ,  10_{143} ,  10_{147}  
  $ 
\end{itemize}

In the next subsection we will observe that 17 of these are concordant to lower genus
knots with $g_c$ known.  We have been unable to resolve the cases of $8_{18}$ and $9_{40}$.

\subsection{Concordances to lower genus knots}

\medskip
{\bf Concordant to trefoil}\qua The following 8 knots are concordant to the trefoil,
$3_1$:
$8_{10}$, $8_{11}$, $10_{40}$, $10_{59}$, $10_{103}$, $10_{106}$, $10_{143}$, $10_{147}$. 

 As noted in
\cite{ta}, in \cite{co} a typographical error leads to the knot $10_{143}$ failing to be
on this list, and $10_{65}$ is placed on the list accidentally.  The knots $10_{103}$,
$10_{106}$ fail to be identified in \cite{co}.  We will describe these concordances below.
 For all of these knots, $ g_c(K) = 1 = g_4(K)$.

\medskip
\noindent{\bf Concordant to the figure eight}\qua  The following 2  knots are concordant to
the figure eight knot, $4_1$:  $9_{24}$,  $9_{37}$.  Both have $g_c(K) = 1 = g_4(K)$.
\medskip

\noindent{\bf Concordant to $5_1$}\qua  The following 2 knots are concordant to the knot
$5_1$:   $10_{21}$,  $10_{62}$. Both have $g_c(K) = 1 = g_4(K)$.
\medskip

\noindent{\bf Concordant to $5_2$}\qua  The following 4 knots are concordant to the knot
$5_2$:        $10_{65}$, $10_{67}$,  $10_{74}$,  $10_{77}$.  The knot $10_{67}$ fails to
be on the list in \cite{co}. Its concordance is described below.  All have
$g_c(K) = 1 = g_4(K)$.

\medskip

\noindent{\bf Concordant to $3_1 \# 3_1$}\qua  The following 1 knot is  concordant to the
connected sum of the trefoil with itself, $3_1 \# 3_1$: $10_{98}$.   It satisfies
$g_c(10_{98}) = 2 = g_4(10_{98})$.

\subsection{Building the concordances for $10_{67}, 10_{103}$ and $10_{106}$}
In this final subsection we will describe the concordances that were not included in
Conway's list.  In Figure 3 we have illustrated the knot $10_{67}$ along
with a band.  Performing the band move  $10_{67}$ along that band results in a split link
of two components: an unknot and the knot $5_2$.  This gives the desired concordance from
$10_{67}$ and $5_2$.

\cfigure{1.75in}{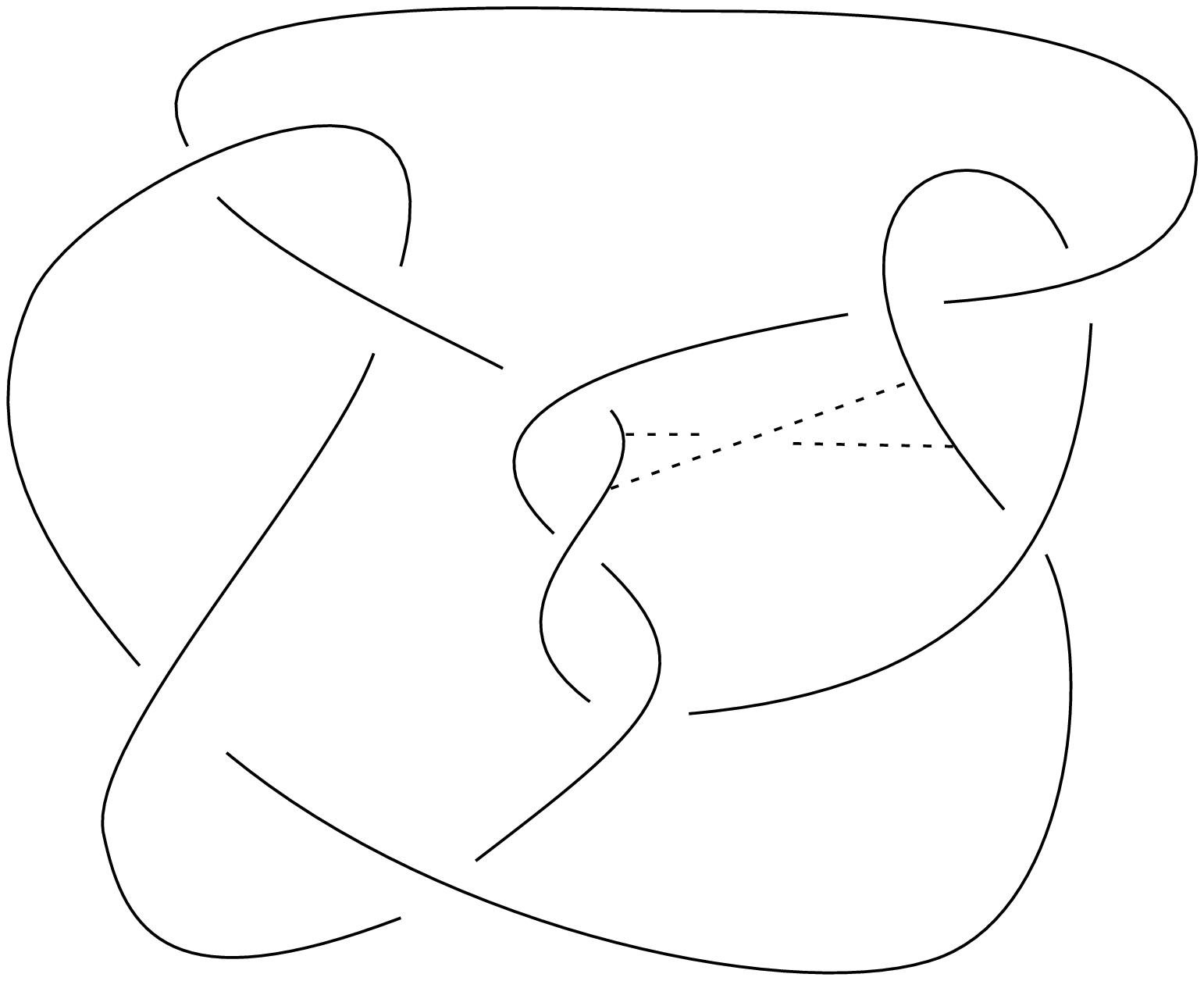}{$10_{67}$  }

In Figure 4 we have illustrated the knots $10_{103} \# {-}3_1$ and  $10_{106}
\# 3_1$  along with a band in each figure.  Performing the band move along each band
yields a split link of two unknotted components.  Hence, both $10_{103} \# -3_1$ and 
$10_{106} \# 3_1$ are slice, so $10_{103}$  and  $10_{106} $ are each concordant to trefoil
knots.

\cfigure{3in}{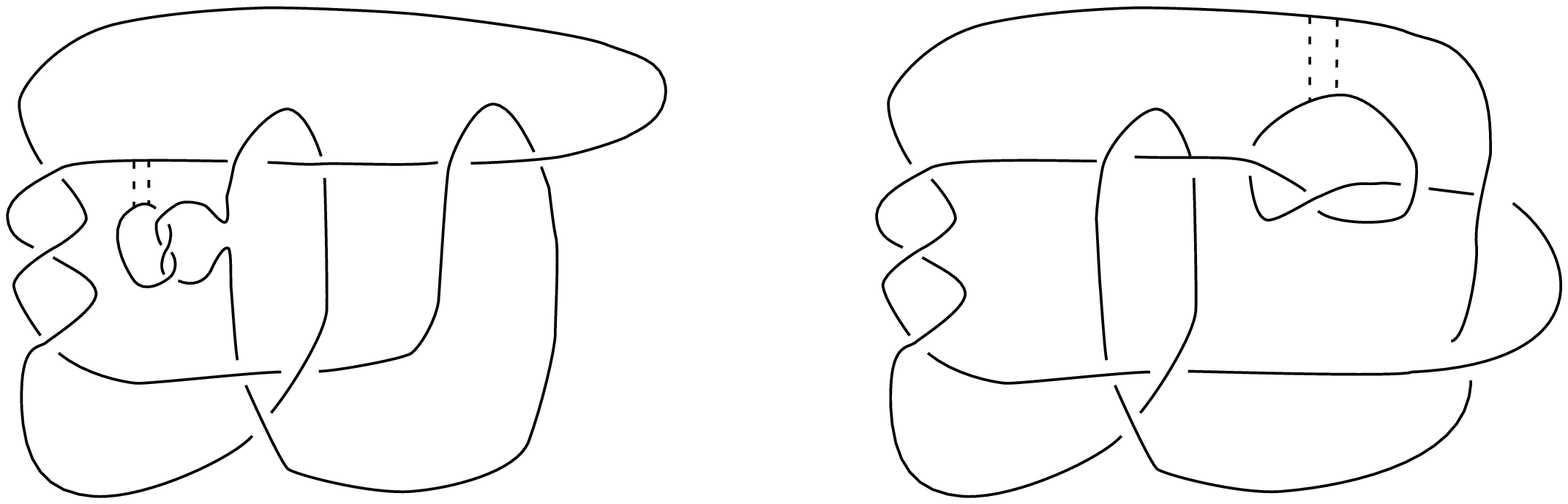}{$10_{103} \# 3_1$\qquad\ \ \ \ \
and\qquad \ \ \ \ \ $10_{106} \# 3_1$\qquad\qquad\quad \hbox{}}

\Addresses\recd

\end{document}

%% file: agtout.tex
\def\ifplaintex{\expandafter\ifx\csname documentclass\endcsname\relax}

\def\gtp{{\mathsurround=0pt\it $\cal G\mskip-2mu$eometry \&\ 
$\cal T\!\!$opology $\cal P\!$ublications}}  

\def\recd{{\small Received:\qua\receiveddate\ifx\reviseddate\relax
\else\qquad Revised:\qua\reviseddate\fi\par}} 

\def\lognumber#1{\def\thelognumber{#1}}
\def\volumenumber#1{\def\thevolumenumber{#1}}
\def\volumeyear#1{\def\thevolumeyear{#1}}
\def\papernumber#1{\def\thepapernumber{#1}}
\def\pagenumbers#1#2{\def\startpage{#1}\def\finishpage{#2}}
\def\published#1{\def\publishdate{#1}}

\def\received#1{\def\receiveddate{#1}}
\def\revised#1{\def\reviseddate{#1}}
\def\accepted#1{\def\accepteddate{#1}}

\long\def\asciiabstract#1{\long\def\theasciiabstract{#1}}

\let\\\par\let\thelognumber\relax\let\thevolumenumber\relax
\let\thepapernumber\relax\let\thevolumeyear\relax\let\startpage\relax
\let\finishpage\relax\let\publishdate\relax\let\receiveddate\relax
\let\reviseddate\relax\let\accepteddate\relax\let\theasciititle\relax
\let\theasciiauthors\relax
\let\theasciiabstract\relax

\let\theasciiemail\relax

\ifplaintex
\font\logobig=cmssbx10 scaled 3836
\font\logomed=cmssbx10 scaled 2557
\else
\font\logobig=cmssbx10 scaled 4200
\font\logomed=cmssbx10 scaled 2800
\fi

\long\def\makeagttitle{   
\count0=\startpage
\agt\hfill      
\hbox to 45truept{\vbox to 0pt{\vglue -13truept{\logomed A\kern -.37em{\logobig 
T}\kern -.38em G}\vss}\hss}
\break
{\small Volume \thevolumenumber\ (\thevolumeyear)
\startpage--\finishpage\nl
Published: \publishdate}

\vglue .25truein

{\parskip=0pt\leftskip 0pt plus
1fil\def\\{\par\smallskip}{\Large\bf\thetitle}\par\medskip} \vglue
0.05truein

{\parskip=0pt\leftskip 0pt plus 1fil\def\\{\par}{\sc\theauthors}
\par\medskip}%
 
\vglue 0.03truein 

{\small\leftskip 25truept\rightskip 25truept{\bf Abstract}\stdspace\theabstract

{\bf AMS Classification}\stdspace\theprimaryclass
\ifx\thesecondaryclass\relax\else; \thesecondaryclass\fi\par
{\bf Keywords}\stdspace \thekeywords\par}\vglue 7truept

}   

\ifplaintex
\font\phead=cmsl9 scaled 950
\font\pnum=cmbx10 scaled 913
\font\pfoot=cmsl9 scaled 950
\headline{\vbox to 0pt{\vskip -4.5mm\line{\small\phead\ifnum
\count0=\startpage ISSN 1472-2739 (on-line) 1472-2747 (printed)
\hfill {\pnum\folio}\else\ifodd\count0\def\\{ }%
\ifx\theshorttitle\relax\thetitle\else\theshorttitle\fi\hfill{\pnum\folio}
\else\def\\{ and }{\pnum\folio}\hfill\ifx\theshortauthors\relax\theauthors
\else\theshortauthors\fi\fi\fi}\vss}}
\footline{\vbox to 0pt{\vglue 0mm\line{\small\pfoot\ifnum\count0=\startpage
\copyright\ \gtp\hfill\else
\agt, Volume \thevolumenumber\ (\thevolumeyear)\hfill\fi}\vss}}
\else
\font=cmsl9 scaled 1050
\font\lnum=cmbx10 
\font=cmsl9 scaled 1050
\makeatletter
\def\@oddhead{{\small\lhead\ifnum\count0=\startpage ISSN 1472-2739 
(on-line) 1472-2747 (printed)\hfill {\lnum\number\count0}\else\ifodd\count0
\def\\{ }\ifx\theshorttitle\relax \thetitle \else\theshorttitle\fi\hfill
{\lnum\number\count0}\else\def\\{ and }{\lnum\number\count0}
\hfill\ifx\theshortauthors\relax 
\theauthors\else\theshortauthors\fi\fi\fi}}\def\@evenhead{\@oddhead}
\def\@oddfoot{\small\lfoot\ifnum\count0=\startpage\copyright\ \gtp\hfill\else
\agt, Volume \thevolumenumber\ (\thevolumeyear)\hfill\fi}
\def\@evenfoot{\@oddfoot}
\makeatother
\fi
\let\maketitlepage\makeagttitle
\let\makeshorttitle\maketitlepage
\let\maketitle\maketitlepage

\newwrite\gtoutfile
\long\gdef\makeheadfile{  
{\def\\{, }\def\s{ }
\immediate\openout\gtoutfile head.xxx
\immediate\write\gtoutfile{Proxy-for: \ifx\theasciiauthors\relax
\theauthors\else\theasciiauthors\fi\s<\ifx\theasciiemail\relax\theemail\else\theasciiemail\fi>}
\immediate\write\gtoutfile{\noexpand\\}
\immediate\write\gtoutfile{Authors: \ifx\theasciiauthors\relax
\theauthors\else\theasciiauthors\fi}
{\def\\{ }\immediate\write\gtoutfile{Title: \ifx\theasciititle\relax
\thetitle\else\theasciititle\fi}}
\immediate\write\gtoutfile{Subj-class: GT or SG, GR etc}
\immediate\write\gtoutfile{MSC-class: \theprimaryclass\ifx\thesecondaryclass\relax\else, \thesecondaryclass\fi}
\immediate\write\gtoutfile{Journal-ref: Algebr. Geom. Topol. \thevolumenumber\s
(\thevolumeyear) \startpage-\finishpage}
\immediate\write\gtoutfile{Comments: Published by Algebraic and
Geometric Topology at}
\immediate\write\gtoutfile{\s\s\s  http://www.maths.warwick.ac.uk/agt/AGTVol\thevolumenumber/agt-\thevolumenumber-\thepapernumber.abs.html}
\immediate\write\gtoutfile{\noexpand\\}
\immediate\write\gtoutfile{}
\ifx\theasciiabstract\relax
\immediate\write\gtoutfile{\theabstract}\else
\immediate\write\gtoutfile{\theasciiabstract}\fi
\immediate\write\gtoutfile{}
\immediate\write\gtoutfile{\noexpand\\}
\immediate\write\gtoutfile{}
\immediate\closeout\gtoutfile}}  

\def\maketitlepage{\makeagttitle\makeheadfile}
\let\makeshorttitle\maketitlepage
\let\maketitle\maketitlepage